\documentclass[11pt]{amsart}
\usepackage{psfrag}
\usepackage[dvips]{graphicx}

\newtheorem{theorem}{Theorem}
\newtheorem{lemma}[theorem]{Lemma}
\theoremstyle{definition}

\newtheorem{remark}[theorem]{Remark}
\numberwithin{equation}{section}

\theoremstyle{corollary}
\newtheorem{example}{Example}
\numberwithin{equation}{section}
\theoremstyle{example}

\theoremstyle{proposition}
\newtheorem{proposition}[theorem]{Proposition}
\newfont{\EUL}{eufm10 scaled 1000}

\newcommand\R{\mathbb{R}}

\newcommand\C{\mathbb{C}}

\renewcommand\P{\mathbb{P}}

\newcommand\g{\mbox{\EUL g}}

\newcommand\SO{{\rm SO}}

\newcommand\SU{{\rm SU}}
\newcommand\U{{\rm U}}

\begin{document}
%
%
\title{On deformations of Hamiltonian actions}
\author{Lucio Bedulli and Anna Gori}
\address{Dipartimento di Matematica \lq U.Dini\rq\\ Viale Morgagni 67/A\\
50100 Firenze\\Italy} \email{bedulli@math.unifi.it}
\address{Dipartimento di Matematica e Appl.\ per l'Architettura\\
Piazza Ghiberti 27\\50100 Firenze\\Italy}
\email{gori@math.unifi.it}
\begin{abstract}
In this paper we generalize to coisotropic actions of compact Lie groups a theorem of
Guillemin on deformations of Hamiltonian structures on compact symplectic manifolds.
We  show how one can reconstruct from
the moment polytope the symplectic form on the manifold.
\end{abstract}
\maketitle
\section{Introduction}
Let $G$ be a compact connected Lie group acting in a Hamiltonian fashion on a $2n$-dimensional compact
symplectic manifold $M$. We will denote by $\mu$ the corresponding moment map from $M$ to the dual of the Lie 
algebra of $G$. Throughout the following we will indicate  Lie groups
and their Lie algebras  with capital and gothic letters respectively. \\
Fix a maximal torus $T$ in $G$ and denote by $\frak{t}^*$ the dual of its Lie algebra.\\
A well known result of Kirwan \cite{Kir2} states that  the moment
map image of $M$ meets the closure of a Weyl chamber $\mathfrak{t}^*_+$ in a
convex polytope $\Delta(M)$. Delzant
\cite{Del2} conjectured that $\Delta(M)$, together with some
additional invariants of the action, determines the manifold up
to $G$-equivariant symplectomorphisms (this has been solved in
special cases, see e.g. \cite{Del1}, \cite{Del2}, \cite{Ig},
\cite{Wood} and \cite{Chiang}). We deal with a local version of
this conjecture. We consider pairs $(\omega,\varphi)$ where $\omega$ 
denotes the symplectic form on $M$
and $\varphi$ the Hamiltonian differentiable action of  a compact  Lie group $G$ on $M$. 
We study smooth deformations of Hamiltonian pairs on a fixed manifold
proving the following
\begin{theorem}\label{main}
Let $(M,\omega)$ be a compact symplectic manifold acted on coisotropically and effectively by a compact  Lie group $G$ in a Hamiltonian
fashion.  The moduli space of Hamiltonian $G$-structures $(\omega,\varphi)$ on $M$ is a smooth manifold whose dimension $k$ is equal to the second Betti
 number of $M$. Moreover there are $2k$ points, not necessarily distinct, $q_1,q_2,\cdots,q_k,q_1',q_2',\cdots,q_k'$ in the moment
 map image
such that, as one varies the pair $(\omega,\varphi)$, the distances between $q_i$ and $q'_i$ are coordinates of this moduli space.
\end{theorem}
\noindent Guillemin in \cite{Gui} solves the problem in the abelian case under the following additional
assumptions
\begin{enumerate}
\item A torus $T$ acts locally freely on an open dense subset of $M$;
\item the set of fixed points,
$M^T$, is finite;
\item the restriction of the moment map on the set $M^T$ is injective.
\end{enumerate}
\noindent In the same paper the author remarks also that one can
treat the non-abelian case replacing the images of $M^T$ with the set of vertices of the Kirwan's polytope  under some
additional hypotheses on the action.  Actually, in our proof we show that the relevant points are exactly the images of points fixed by
a previously chosen maximal torus of $G$.

\vspace{0.3cm}

The paper is organized as follows: we firstly prove a version of Guillemin's theorem in case  the moment map restricted to ${M^T}$ {\em is not injective}. Then,
with an elementary argument, we show how to recover the image of the moment map
for the action of the maximal torus from the Kirwan's polytope (Lemma \ref{hull}). We finally show that if the $G$-action is coisotropic
the problem boils down to the abelian situation treated in
\cite{Gui} in the case the map $\mu_{|_{M^T}}$ is not injective.\\
\noindent {\bf Acknowledgements.} The authors would like to thank
Prof. V. Guillemin and Prof. Elisa Prato for their advices.
\section{The abelian case}
\label{sez2}
\begin{figure}[h]
\centering
\includegraphics[height=5.5cm]{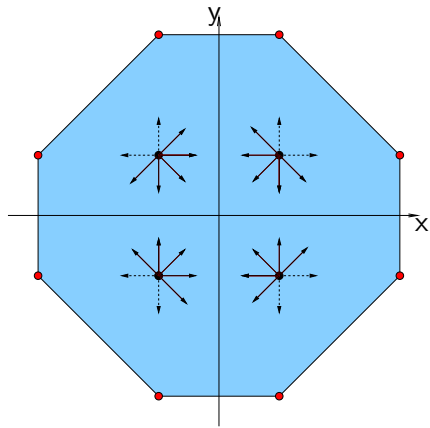} \caption{Moment polytope of a
$T^2$-action on the product of two coadjoint orbits of $\SO(5)$.}
\label{controesemp}
\end{figure}

Throughout this section a Hamiltonian action of a torus $T$ on a
compact symplectic manifold $(M,\omega)$ is fixed. Denote by
$\mu:M\to \mathfrak{t}^*$ the corresponding moment map.\\ For the
reader's convenience we here recall some notations from \cite{Gui}
adapted to our purposes.  Let $p$  be a $T$-fixed point, and let
$$\Sigma_p=\{\alpha_1^p,\alpha_2^p,\ldots,\alpha_n^p\}\subset\mathfrak{t}^*$$
be the set of weights of the isotropy representation on $T_p M$.
From now on we assume that {\em $M^T$ is finite}. Choose a topological generator of $T$, say
$\xi\in\mathfrak{t}$ 
 such that $\alpha_i(\xi)\neq 0$ for all $\alpha_i\in \cup \Sigma_p$ where $p$ ranges through $M^T$.
Denote by $\sigma_p$ the cardinality of the set $\{\alpha \in \Sigma_p: \alpha(\xi)<0\}$.
Consider now the real valued function $\mu_{\xi}: x\mapsto \mu(x)(\xi)$. The defining properties of $\mu$ implies that the
critical points of $\mu_{\xi}$ are the fixed points of $T$. Moreover  $\mu_{\xi}$ is a  Morse function (see e.g. 
\cite{Atiyah} and \cite{Fra}) and the index  of a critical point $p\in M^T$ equals $2\sigma_p$.\\
The fact that there are no critical points of odd index  implies that
$\mu_{\xi}$ is  perfect,  hence the even Betti numbers  $b_{2k}$ equals the number of points in $M^T$ with
$\sigma_p=k$ and the odd ones vanish.\\
\noindent The moment map image of $M$ is convex \cite{Atiyah}, and
the directions of the edges starting from a vertex $q\in\mu(M^T)$
are given by the weights of the representations at points $p_i$ in
the fiber of $q$. More precisely, let $q$ be the image of a
$T$-fixed point $p$ in $M$, via the moment map. On every direction
$\alpha\in \Sigma_p$ starting from $q$ there exists another point
$q'\in \mu(M^T)$ \cite{Gui}. It is natural to investigate whether
one can or not ``determine'' the even Betti numbers of $M$ from
the image of the moment map.\\ When the restriction of the
moment map to $M^T$ is not injective {\em this is not always
possible}; if two points $p,p'$ have the same image $q$ via the
moment map they describe the same cones
$C_p=cone(\alpha_1^p,\alpha_2^p,\ldots,\alpha_n^p)$ and $
C_{p'}=cone(\alpha_1^{p'},\alpha_2^{p'},\ldots,\alpha_n^{p'})$
(this can be seen combining Example 4.21 and  Theorem 6.5 in
\cite{Sja} since the fibers of the moment map are connected)
nevertheless they can have $\sigma_p\neq \sigma_{p'}$, for any choice of the topological generator of $T$, as Example
\ref{controesempio} shows.
\begin{example}\label{controesempio}
{\em Let $T$ be a maximal torus in $\SO(5)$. The usual choice for a basis for $\mathfrak{t}$ gives isomorphisms $\mathfrak{t}\cong \R^2$ and $$\mathfrak{t}^*_+\cong \{(x,y)\in \R^2\::0\leq y\leq x \}.$$
Let $\Theta_{\gamma}$ be the coadjoint orbit through $(\gamma,\gamma)$ 
and  $\Theta_{\delta}$ be  the coadjoint orbit through $(\delta,0)$ where $\gamma$ and $\delta$ are positive real numbers.
Consider the diagonal action of $T$ on the product $M=\Theta_{\gamma}\times \Theta_{\delta}$.
If $\gamma=\delta$ the image of the moment map is given in Fig. \ref{controesemp}, and the dots are the image of the $T$-fixed points.
The interior dots represent points that have two inverse images. One can determine the weights of the representations
on the tangent spaces at two different points in the same fiber.
As one can check from the picture, they differ for one weight, hence for any choice
of $\xi$ there is always a pair of $T$-fixed points with the same image and different  $\sigma$.
In the figure the filled arrows denote the common weights, while the
dashed ones come from different pre-images of the same point.}
\end{example}
\psfrag{x}{$\sf x$}
\psfrag{y}{$\sf y$}
The previous example clearly shows that one cannot determine the
Betti numbers of $M$ starting from $\mu_T(M)$ without coming back
to the manifold and counting the contribution of every point in
$M^T$. Nevertheless, with the same arguments used in \cite{Gui}
(pag. 230 Lemma $2$), it is always possible to associate to each point in
$p\in M^T$ of index $k$ an element $c_p \in H_{2k}(M)$ and show
that the set of $\{c_p,\;\;\sigma_p=k\}$ yields a {\em basis} of
$H_{2k}(M)$.
\\Now consider, in particular, the set of points $p_i\in M^T$ with
$\sigma_{p_i}=1$; the images of these points can be distinct or
not. Note that for each $p_i$ in this set there exists only one
$\alpha_j^{p_i}\in \Sigma_{p_i}$ with $\alpha_j^ {p_i} (\xi)<0.$
Moreover on the ray
$$q+\alpha_j^{p_i}t\;\;\;0<t<\infty$$ there is at least another
point $q'\in\mu(M^T)$. More precisely, we have
\begin{lemma} \label{distanze} Let $t_0=[\omega](c_p)$. Then the point $q'=q+t_0
\alpha_j^{p}$, where $p$ has $\sigma_p=1$, belongs to $\mu(M^T)$.
\end{lemma}
\noindent This can be proven following the same ideas of
\cite{Gui}. Thus, in particular, if $p_1,p_2\ldots, p_k$ are
points of index $1$, the distances between their (not necessarily
distinct) images, via the moment map, $q_1,q_2\ldots, q_k$ and
$q'_1,q'_2\ldots, q'_k$ determine the cohomology class of $\omega$
in $H^2(M)$.\\
Now we recall some definitions and facts on deformations. By a
differentiable {\em deformation} of a differentiable action
$\varphi$ of a group $G$ on a manifold $M$ we mean a one-parameter
family $\varphi_t$ ($t\in I=[0,1]$) of differentiable actions such
that $\varphi_o=\varphi$ and the map
$(g,p,t)\mapsto\varphi_t(g,p)$ of $G\times M\times I$ into $M$ is
differentiable. Recall that a deformation of $M$ is a
one-parameter family $\psi_t$ of diffeomorphisms of $M$ such that
$\psi_0$ is the identity and the map  $(p,t)\to \psi_t(p)$ is
differentiable. A deformation action $\psi_t$ of a manifold $M$ under
the action $\varphi$ of a group $G$ induces a differentiable
deformation $\varphi_t$ defined by 
\[
\varphi_t(g,p):=\psi_t(\varphi(g,\psi_t^{-1}(p))).
\]
Such a deformation $\varphi$ is called {\em trivial}.  Palais and Stewart showed that
\begin{proposition}\cite{PS} If
the group $G$ and the manifold $M$ are assumed to be compact then
any differentiable deformation of a given $G$-action $\varphi$ on
$M$ is differentiably trivial.
\end{proposition}
Hence we can assume that the $G$-action is fixed. Moreover, if we
``move'' the symplectic form in its cohomology class, using
a $G$-equivariant version of Moser's Theorem, we get that the moment map image 
does not change.
Henceforth the only way to deform the Hamiltonian pair $(\omega,\varphi)$ is
to deform the class $[\omega]$.\\ 
After a deformation of the class of $\omega$, the
directions of the edges of $\mu_T(M)$ are preserved,  while the distances of certain
points give us a ``measure'' of the deformation. Indeed, if we put
on the manifold $M$ a $T$-invariant metric $g$ and consider the
map that associates to each class $c$ in $H^2(M,\R)$ the $2$-form
$\omega_c=\omega+h^{g}_{c-c_0}$, where $c_o$ is the cohomology
class of $\omega$ and $h^{g}_{c-c_0}$ the unique harmonic
representative of $c-c_o$, we get a $T$-invariant closed $2$-form,
$\omega_c$, which is symplectic if $c$ is close to $[\omega]$.
Moreover, since $H^1(M,\R)=0$, the action $\varphi$
is Hamiltonian on $(M,\omega_c)$. Therefore {\em the set of pairs
$(\omega_c,\varphi)$, with $c$ close to $[\omega]$, is an open
subset of the moduli space of $G$-Hamiltonian actions on $M$ whose
coordinates are thus given by the distances of $q_1,q_2\ldots,
q_k$ and $q'_1,q'_2\ldots, q'_k$}, according to Lemma
\ref{distanze} and the subsequent remarks.
\section{The non abelian case: proof of Theorem 1}
In this section we consider a compact, non necessarily abelian, Lie group $G$
acting in a Hamiltonian fashion on a compact symplectic manifold $(M,\omega)$. Fix once and for all a maximal torus $T$ in $G$.
We will denote by $\mu$ and $\mu_T$ the corresponding moment maps for the $G$- and the $T$-action respectively.
The following lemma shows how to recover the $\mu_T$ image from the Kirwan polytope $\Delta(M)$,
more precisely from $\mu(M) \cap \frak{t}^*$ (which is obtained from
$\Delta(M)$ by reflections through the walls of the Weyl chambers). We
state and prove it since we could not find it in the literature,
although it is probably well-known.
\begin{lemma}
\label{hull}
The set $\mu_T(M)$ coincides with the convex hull of $\mu(M) \cap \frak{t}^*$.
\end{lemma}
\begin{proof}
Firstly recall that the image of the moment map for the $T$-action is the convex hull
of the finite set $\mu_T(M^T)$ \cite{Atiyah}. Note that $\mu$ and $\mu_T$ coincide on $M^T$ indeed the restriction
of the moment map $\mu$
to the fixed point set of a closed subgroup $H$ of $G$, takes values
in the dual of the Lie algebra of the centralizer of $H$ in $G$, $\frak{c_G(h)}^*$, hence
for $H=T$ the image is contained in $\frak{t}^*$.
Therefore $\mu_T(M)={\rm{conv}}(\mu(M^T)) \subset {\rm{conv}}(\mu(M)\cap \frak{t}^*)$.
The other inclusion follows from the fact that
$\mu_T$ is the
composition of $\mu$ with the projection $\frak{g}^* \to \frak{t}^*$
induced by the inclusion $\frak{t} \hookrightarrow \frak{g}$.
\end{proof}
The previous lemma is still true if the moment map $\mu$ is assumed to be proper, and the manifold is not necessarily compact.\\
If we assume that the $T$-fixed point set is 
finite, and that $T$ acts on an open dense subset $M_0 \subseteq M$ locally freely, 
we can apply the results of section \ref{sez2} to the non-abelian case. A large class of actions for which $M^T$ is finite is provided by {\em coisotropic} ones.
Recall that a Hamiltonian $G$-action is called {\em coisotropic}  if the principal
$G$-orbits $G \cdot p$ are coisotropic with respect to $\omega$ i.e.
$(T_p \;G\cdot p)^\omega\subseteq T_p \;G\cdot p$.
Before proving our main theorem we here state a theorem that gives a
characterization of coisotropic actions on non necessarely compact symplectic manifolds. 
We will use it only in the compact setting, however we give here the 
proof in full generality since it might have an autonomous interest (see also \cite{wood} in the compact setting).
\begin{theorem} 
{\label{class}}The following conditions for a Hamiltonian $G$-action on a connected symplectic manifold $M$ 
with connected moment map fibers are equivalent :
\begin{itemize}
\item [(i)] The $G$-action is coisotropic.
\item [(ii)] The cohomogeneity of the $G$ action is equal to the
difference between the rank of $G$ and the rank of a principal
isotropy subgroup of $G$.
\item [(iii)] The set of $G$-invariant smooth functions $C^\infty(M)^G$ is an abelian Poisson algebra.
\item[(iv)] The moment map $\mu$ separates $G$-orbits.
\end{itemize}
\end{theorem}
We here  enclose condition (ii) as it can be used  in proving
that
Example \ref{stellette} in Section $4$ is coisotropic.
\begin{proof} We get our claim, combining a result of \cite{HW} and a result of \cite{ACG}.\\
The first two conditions are equivalent thanks to Theorem 3.3 in \cite{HW}. The equivalence of (i) and (iii) is straightforward.
We want to show that if (iii) holds then all the 
reduced spaces $M_q=\mu^{-1}(q)/G_q=\mu^{-1}(G \cdot q)/G$ are points. In general $M_q$ is not necessarely 
smooth, however  Theorem $2$ in \cite{ACG} says that, for arbitrary $q$ the algebra 
$$C^\infty(M_q):=C^\infty(M)^G/{\mathcal{I}}_q^G$$
where ${\mathcal{I}}_q^G$ is the ideal of $G$-invariant functions vanishing on $\mu^{-1}(q)$, is a non-degenerate Poisson algebra. 
That is the Poisson bracket vanishes only 
on constant functions on $M_q$. If $C^\infty(M)^G$ is abelian, then $C^\infty(M)^G/{\mathcal{I}} _q^G$ is 
abelian and hence $M_q$ is a point, whenever it is connected (this naturally holds when the fibers of $\mu$ are connected). Finally, (i) follows from (iv) immediately because (i) means 
that the moment map fibers are generically tangent to $G$-orbits.  
\end{proof}
We are now ready to prove Theorem \ref{main}
\begin{proof}[Proof of Theorem 1]We firstly show that {\em if the action is coisotropic then the $T$-fixed
points are isolated}. Since $\mu(M^T)$ is finite and $\mu$
separates  orbits (iv), the set $\mu^{-1}(\mu(M^T))$ is
contained in the union of finitely many $G$-orbits. The claim
follows from the fact that the number of $T$-fixed points in a
$G$-orbit $\mathcal O$ is given by its Euler characteristic
$\chi(\mathcal O)$.\\ \noindent Now we only need the $T$-action to
be locally-free on an open dense subset of $M$. In fact this is
the case when the $G$-action is effective. Indeed, if $M_o
\subseteq M$ is the dense open subset of $G$-prinicipal points,
and $L$ is a principal isotropy subgroup of $T$ acting on $M_o$,
we have that $L$ acts trivially both on the tangent space to the
orbit and on the slice,
therefore it is contained in the kernel of the action.\\
Finally we can apply the non-injective version of Guillemin's
theorem, in order to obtain the coordinates of the moduli space of
$T$-invariant Hamiltonian structures on $M$. Note that, since we start
from a $G$-invariant symplectic form $\omega$, we can choose a
$G$-invariant metric $g$ on $M$, compatible with $\omega$, thus
the $g$-harmonic forms $h^g_{c-c_o}$ will be $G$-invariant, and
the claim  follows.
\end{proof}
In the following example in which the action is not coisotropic, and the $T$-fixed point set 
is not finite, Theorem $1$ does not hold.
\begin{example}{\em
Consider the manifold $M$ given by the product of a full flag $G/T$  and a manifold $N$ acted on transitively 
on the first factor and trivially on the second one by $G$. 
The moment polytope for this action is a point. The convex hull of the
points obtained by reflections through the Weyl walls gives us
$\mu_T(M)$. The ``deformations parameters'' are at most
$l=\mbox{rank}(G)$, while the second Betti number can be larger, whenever $b_2(N)>l.$}
\end{example}
\section{examples}
We have emphasized that the relevant points in the moment polytope in order to study deformations of Hamiltonian
 structures are the images of the $T$-fixed points. It would be interesting to distinguish from the combinatoric of the moment image 
the points in $\mu(M^T)$. We here collect some facts in this regard.\\
We say that a convex polytope $\Delta(M)\subset \mathfrak{t}^*_+$ is {\em reflective} \cite{Wood} at $q\in \Delta(M)$ if and only if
\begin{enumerate}
\item The set of hyperplanes that intersect $\Delta(M)$ in codimension
 $1$ faces that contain $q$ is invariant under the stabilizer $W_q$ of $q$, where $W$ is the Weyl group of $G$;
 \item any open codimension $1$ face of $\Delta(M)$ containing $q$ in its closure is contained in the open positive
 Weyl chamber in $\mathfrak{t}^*_+.$
 \end{enumerate}
\noindent The first condition is equivalent to requiring that if $H$ is an hyperplane
such that $H\cap \Delta(M)$ is a codimension $1$ 
face containing $q$ and $w\in W_q$ then $H\cap w\Delta(M)$ is a
codimension $1$ face of $w\Delta(M)$.
That is, the codimension $1$ faces continue through the walls.\\\\
$\bf(1)$ \label{riflessivi} Let $q\in\mu(M)$  be a non reflective
point. Then $q$ is the image of a $T$-fixed point.\\
Let $q$ be a non reflective vertex on a wall. Let $\alpha\in \mathfrak{t}^*$,  be a direction, 
starting from $q,$  such that it ``does not continue'' beyond the
wall. 
Let $\mathfrak{k}=\{\eta\in \mathfrak{t},\: \alpha(\eta)=0\}$. 
Since $\alpha$ belongs to $\mathfrak{t}$,  $\mathfrak{k}$ is the lie
algebra of a closed connected  Lie subgroup, $K$, of $T$;
the quotient group $T/K$ is a circle group with Lie algebra
$\mathfrak{t}/\mathfrak{k}$. Let $X$ be the connected component of $M^K$ 
that is sent by $\mu$ onto the line $q+t\alpha$. The group $K$ acts
trivially on $X$ and the $T/K$-action on it is Hamiltonian 
with moment map $\mu_\alpha:X\to\R$ defined by $\mu_\alpha(x)=\langle
\mu(x),\alpha \rangle$. The critical points of
$\mu_\alpha$ are those $p\in X$ where the fundamental field associated 
to $\alpha$ vanishes, i.e. points fixed by the torus $T^\alpha$
generated by $\alpha$.
The function $\mu_{\alpha}$ achieves its maximum and its 
minimum on the extremal points 
of the segment $\{q+t \alpha: t \in [0.t_0]\}$. Therefore $q$ is the image of a point
$p\in M^K$ critical for $\mu_\alpha,$ and  it is fixed by $T$.\\\\
$\bf(2)$ The interior vertices $q$ ,i.e.\ vertices that belongs to the
interior of a Weyl chamber, 
are reflective ($W_q$ is trivial) and images of $T$-fixed points.\\ Since in this case the stabilizer $G_q$ is exactly the maximal torus
$T$, the last statement is a consequence of the necessary condition for a  point $p$ of $M$ to be sent
to a vertex (see pag.81 Theorem $6.7$ \cite{Sja})
\[
\g_q=[\g_q,\g_q]+\g_p,
\] 
from which one gets $G_p=G_q=T$.\\\\
\noindent Furthermore observe that, if the action is coisotropic, the fiber of an interior vertex contains only one $T$-fixed point $p$, 
indeed each $T$-fixed point in the orbit $G\cdot p$  is sent to a different Weyl chamber. Note 
that there can be points $p$ with $G_p=T$ that are sent to points on the boundary of a Weyl chamber (see fig. \ref{AA}).\\
The maximal points in $\mu(M)$, i.e. points whose distance from the origin is maximal, are non reflective henceforth belong to $\mu(M^T)$.
In the K\"ahler setting, this can be seen 
also as a consequence of the fact that if  $p$ is a point in which $\|\mu\|^2$ reaches its maximum,
the $G$-orbit through $p$ is complex \cite{GP}. We remark here that if $p$ is a point in which $\|\mu\|^2$ achieves its minimum 
the $G$-orbit through $p$ should  be symplectic but not complex (see Remark \ref{min} for an example).

\begin{example}{\em  Consider the natural action of $\SU(n+1)$ on the complex $n$-dimensional 
projective space $\P^n$ $$A\mapsto\{[x]\mapsto[Ax]\}.$$  This action preserves the Fubini-Study K\"ahler form $\omega_0$ 
of $\P^n$ and it is obviously Hamiltonian. If we identify the Lie algebra $\mathfrak{su}(n+1)$ with its dual  by means 
of the Killing form, the corresponding moment map (which is unique, due to the semisemplicity of the special 
unitary group  $\SU(n+1)$) can be written as 
follows $$\mu:\P^n\to\mathfrak{su}(n+1)$$ $$[x]\mapsto \frac{1}{2\pi\:i\|x\|^2}\Big({xx^*}- \frac{\mbox{Tr}(xx^*)}{n+1}I_{n+1}\Big)$$
we can obtain a Hamiltonian action of $\SU(n+1)$ on the symplectic product of $\P^n$ with itself in two ways: the natural one
$$A\in\SU(n+1)\mapsto\{([x],[y])\mapsto ([Ax],[Ay])\}$$
or the skew one
$$A\in\SU(n+1)\mapsto\{([x],[y])\mapsto ([Ax],[\bar{A}y])\}.$$
{\bf (A)} We consider, for simplicity, the case $n=2$ and we deal first with the natural action. The case $n>2$ can be treated with the same arguments.\\
The principal orbits have codimension $1$ in $M=\P^2\times\P^2$ 
(they are diffeomorphic to $\SU(3)/U(1)$), hence the action is naturally 
coisotropic. Note that the action of $\SU(n+1)$ is not locally free, but since it is effective, 
the action of a maximal torus is locally free on an open dense subset of $M.$ 
Since $b_2(\P^2\times\P^2)=2$ we can find a $2$-parameter family of $\SU(3)$-invariant symplectic forms on $M$.
 Denote by $\omega_{t,s}$ the $2$-form $t\omega_0\oplus s\omega_0$. For $(t,s)$
close to $(1,1)$ this is symplectic on $\P^2\times\P^2$. The corresponding moment 
map $\phi^{t,s}:\P^2\times\P^2\to\mathfrak{su}(3)$ is given by $$\phi^{t,s}(x,y)=t\mu(x)+s\mu(y).$$
To draw the picture of Kirwan's
 polytope it suffices to find the image of two points belonging to the two singular orbits, and then consider 
the maximal torus of $\SU(3)$ such that the dual of its Lie algebra contains the two images. Take $p_1=\{[1:0:0],[1:0:0]\}$ 
in the complex orbit $\mathcal {O}_1=\{([x],[x]), x\in \C^3\}$ (which is biholomorphic to $\P^2$) and $p_2=\{[1:0:0],[0:1:0]\}$ in 
$\mathcal {O}_2=\{([x],[y]),\::\ ,<x,y>=0 \:x,y\in \C^3\}$ (here $<\cdot,\cdot>$ denotes the standard Hermitian product in $\C^{3}$).
 Then $$\phi^{t,s}(p_1)=\frac{(t+s)}{6\pi\,i}\left(
\begin{array}{ccc}
2&0&0  \\
0&-1&0 \\
0&0&-1 \\
\end{array}
\right)$$
and $$\phi^{t,s}(p_2)=\frac{1}{6\pi\,i}\left(
\begin{array}{ccc}
2t-s&0&0  \\
0&2s-t&0 \\
0&0&-t-s \\
\end{array}
\right)
$$
Note that for $t=s$  the image of $p_2$ lies on a Weyl wall.
When we deform the symplectic structure in such a way that $t\neq s$, $q_2=
\phi^{t,s}(p_2)$ lies in the interior of a Weyl chamber; therefore $p_2$ must be fixed by $T$. By reflection through the wall 
we can find another point $q_3=\phi^{t,s}(p_3)$ where $p_3\in M^T$. Note that $\phi^{t,t}(p_2)=\phi^{t,t}(p_3)$, hence when $t=s$ 
the moment map is not injective on $M^T$. For the $T$-generator $\xi$ opportunely chosen $\sigma_{p_2}=\sigma_{p_3}=1$
 and the distances $d(\phi^{t,s}(p_2),\phi^{t,s}(p_1))$ 
and $d(\phi^{t,s}(p_3),\phi^{t,s}(p_1))$ are coordinates of the moduli space of Hamiltonian structures on 
$\P^2\times\P^2$.\\ Note that the orbit through $p_2$ does not inherit the  symplectic structure $\omega$ of $M$; 
indeed suppose by contradiction that $\omega$ is non degenerate on 
$\mathcal O _2=\SU(3)\cdot p_2$ now the restriction of $\phi^{t,t}$ on
$\mathcal O_2$ should be a covering on  $\SU(3)\cdot {\phi^{t,t}(p_2)}$,
hence  $\SU(3)_{p_2}=T$ should be equal to $\SU(3)_{\phi^{t,t}(p_2)}=\U(2)$
 which is not the case. 
\psfrag{fip1}{$\sf \phi_{t,s}(p_1)$}
\psfrag{fip2}{$\sf \phi_{t,s}(p_2)$}
\psfrag{fip3}{$\sf \phi_{t,s}(p_3)$}
\begin{figure}[h]
\centering
\includegraphics[height=5.5cm]{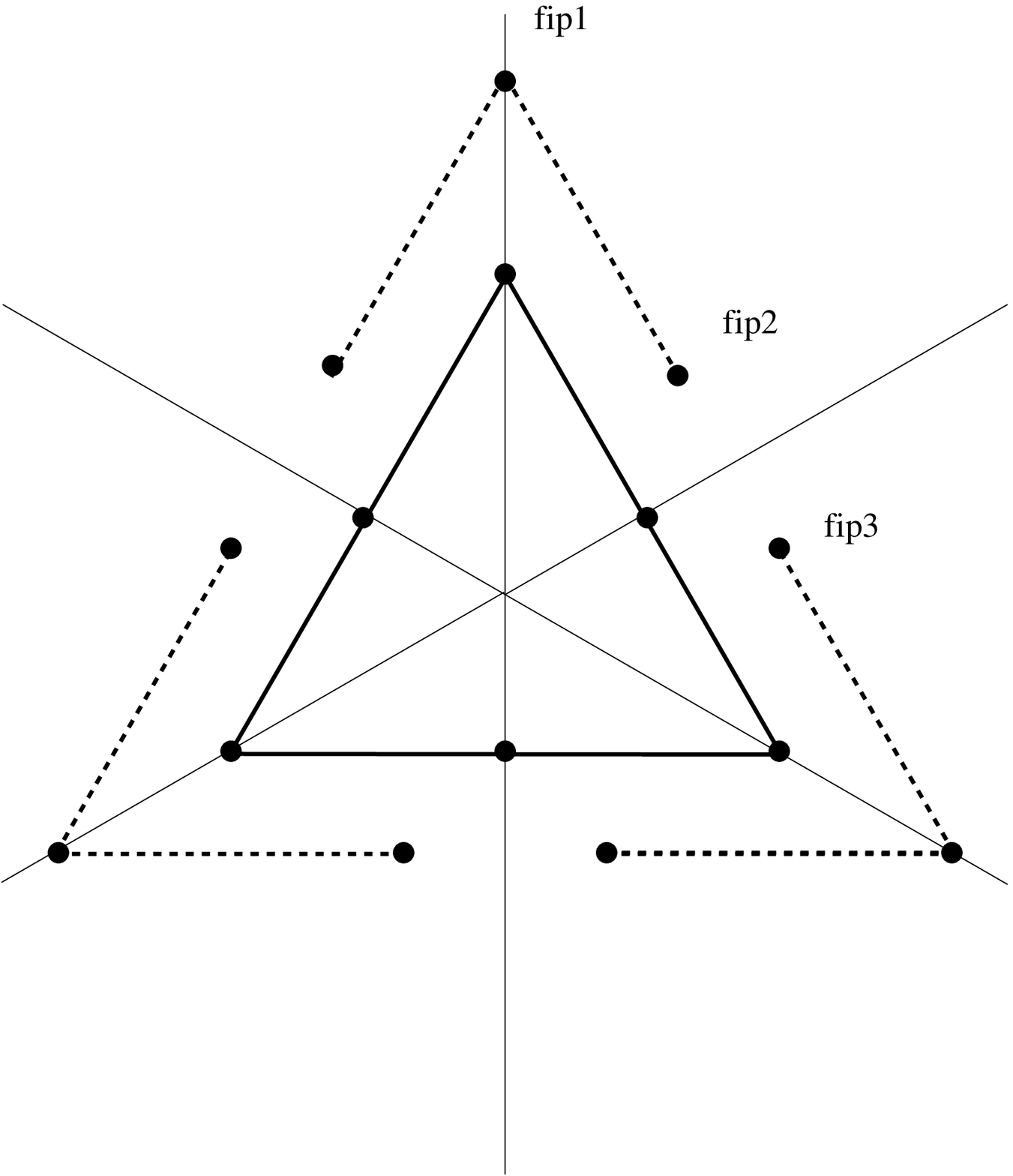} \caption{The moment map image for the  natural action of $\SU(3)$ on $\P^2\times\P^2$:
the continuous lines are the edges when $t=s$, the dashed ones when $t\neq s$.}
\label{AA}
\end{figure}\\ 
{\bf (B)} The picture changes drastically for the skew action of $\SU(3)$ on the same symplectic manifold. 
Once again the principal stabilizer is $1$-dimensional, and the action is coisotropic. The moment map
$\widetilde{\phi}^{t,s}:\P^2\times\P^2\to\mathfrak{su}(3)$ w.r.t.
 the symplectic form  $\omega_{t,s}$ is given by $$\widetilde{\phi}^{t,s}(x,y)=t\mu(x)-s\mu(y).$$ 
 In this case the two singular orbits are  $\mathcal {O}_1=\{([x],[x]), x\in \C^3\}$ and
$\mathcal {O}_2=\{([x],[y]): <x,\bar{y}>=0 \:x,y\in \C^3\}$. We can therefore determine $\Delta(M)$, choosing
 the same points $p_1$ and $p_2$, as before. We see that at least three points $p_2,p_3,p_4$ of $M^T$ are sent by 
$\widetilde{\phi}^{t,t}$ to $0$. When $t\neq s$ the images of these three points move on the walls away from the origin, 
note that for $t=s$ the preimage of $0$ is totally real in $\P^2\times\P^2$, while this is no more true for 
the preimages of the three points on the walls. Since the four  weights of the representation of $T$ on $T_{p_i}M$ $i=2,3,4$
occur in opposite pairs (both $\alpha$ and $-\alpha$ are weights), $\sigma_{p_i}=2$ 
for any choice of the topological generator $\xi$. 
Hence the two relevant vertices (with index $\sigma=1$) come from two distinct points of the complex orbit $\mathcal {O}_2$.
\begin{figure}[h]
\centering
\includegraphics[height=5.5cm]{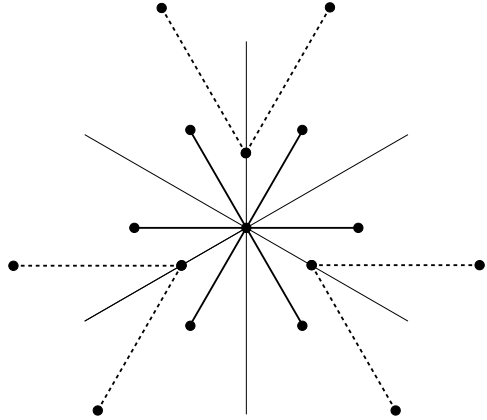} \caption{The moment map images for the skew action of $\SU(3)$ on $\P^2\times\P^2$:
the continuous lines are the edges when $t=s$, the dashed ones when $t\neq s$.}
\label{AAbar}
\end{figure}\\ 
In the previous example, when considering the diagonal action of $\SU(3)$ on the K\"ahler manifold
$(M,\omega)=(\P^2\times\P^2,\omega_o\oplus \omega_o)$, we remarked that  $\mathcal {O}_1$ is 
a codimension $2$ complex orbit (it corresponds to a maximum).
 Then, we can blow up $M$ along  $\mathcal {O}_1$ in a $\SU(3)$-equivariant manner, 
obtaining a K\"ahler manifold  $\widetilde M$ whose second Betti number is $3$.  
Note that the cohomogeneity of the $\SU(3)$-action on $\widetilde M$ is 
unchanged and the action is therefore coisotropic.
Using the technique
of symplectic cuts (see e.g.\ \cite{Wood}, and \cite{Ler}) it is not hard to draw
the moment polytope as in figure \ref{blow}. Note that now the ``degrees 
of freedom''  are $3$ (the ``blown up vertex '' has left the Weyl wall). The 
singular orbits in $\widetilde M$ both have codimension $1$, then a further blow-up will not affect $b_2$.    
\begin{figure}[h]
\centering
\includegraphics[height=5.5cm]{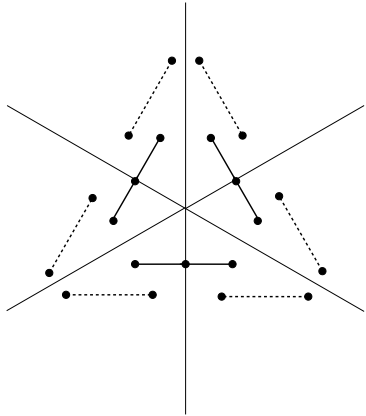} 
\caption{The moment map image for the  natural action of $\SU(3)$ on the blow-up of  $\P^2\times\P^2$ 
along a singular orbit of codimension $2$:
the continuous lines are the edges when $t=s$, the dashed ones when $t\neq s$.}
\label{blow}
\end{figure}\\ 
}
\end{example}
\begin{remark}\label{min}
The $G$-action on the $G$-equivariant blow-up $\widetilde M$ of a K\"ahler manifold $M$ acted on 
coisotropically by a group of isometries $G$, 
along a complex orbit $\mathcal O$, is still coisotropic. \\
 Recall that, in this setting, when $G$ is a
compact Lie subgroup of the full isometry group, if all Borel
subgroups of $G^\C$ act with an open orbit on $M$, then the
$G^\C$-open orbit $\Omega$ is called a {\em spherical homogeneous
space} and $M$ is called a {\em spherical embedding} of $\Omega.$
The $G$-action is coisotropic if and only if the K\"ahler manifold $M$ is projective algebraic,
$G^\mathbb{C}$-almost homogeneous and a spherical embedding of the
open $G^\mathbb{C}$-orbit \cite{HW}.\\
Recall that on the $G$-equivariant blow-up $\widetilde M$ of a K\"ahler manifold $M$ there is a well defined $G$-invariant 
K\"ahler structure (see \cite{Bl}). Denote by $\pi:\widetilde M\to M$ the canonical projection,
by $\widetilde{\mathcal O} = \pi^{-1}(\mathcal{O})$ the exceptional divisor 
and by $\Omega$ the open (hence dense) orbit of a Borel subgroup $B$ of $G^\C$ in $M$.
Since $\mathcal O$ is closed in $M$ and the restriction 
$\pi|_{\widetilde{M} \setminus \widetilde {\mathcal O}}: \: {\widetilde{M} \setminus \widetilde{\mathcal O}} \to M \setminus {\mathcal O}$
is a biholomorphism, the submanifold $\widetilde
\Omega=\pi^{-1}(\Omega)$ is a $B$-orbit open in  
${\widetilde{M} \setminus \widetilde{\mathcal O}}$ 
hence also in $\widetilde{M}$.\\
We conclude this remark pointing out that a symplectic blow-up
$\widetilde M$ of a K\"ahler  manifold $M$ along 
a non complex orbit, does not necessarely admit an invariant compatible K\"ahler structure 
as Theorem $7.3$ \cite{Wood} shows. In particular, one can argue that 
the $G$-orbit through a minimal point of $\|\mu\|^2$ is not complex in this case. 
\end{remark}
\begin{example}
{\em
Now we go back to Example \ref{controesempio}, 
and we consider the action of the semisimple group $\SO(5,\R)$ instead of its maximal torus. 
We refer to \cite{Wood}, in drawing the picture on the left. With the same notations of Example \ref{controesemp}
 we consider the non transversal case ($\gamma=\delta$).
Note that the image in Example \ref{controesemp} can be drawn taking the convex hull of this picture as Lemma \ref{hull} says.
By Kunneth's formula the second Betti number of $M$ is $2$, thus, once again deformations have two degrees of freedom. Note that
 the vertex $q_{t,s}\notin \mu(M^T)$, hence it never leaves the wall (remember that interior vertices come from $M^T$ see fact $\bf(2)$). 
The vertex $r_{t,t}$ is the image of a point in  $M^T$ since is non-reflective, see fact $\bf(1)$. 
 
\psfrag{qts}{$\sf q_{t,s}$}
\psfrag{rts}{$\sf r_{t,s}$}
\psfrag{sts}{$\sf s_{t,s}$}
\psfrag{q11}{$\sf q_{t,t}$}
\psfrag{r11}{$\sf r_{t,t}$}
\psfrag{s11}{$\sf s_{t,t}$}

\begin{figure}[h]
\centering
\includegraphics[height=5.5cm]{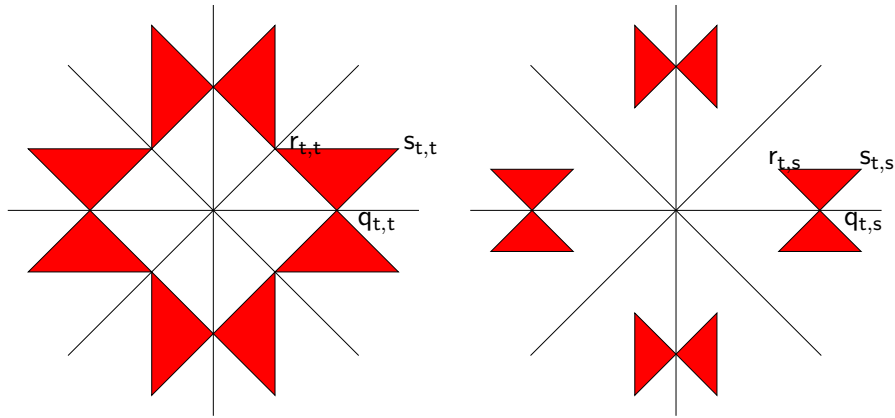} 
\caption{The moment map images for the action of $\SO(5)$ on 
$\frac{\SO(5)}{\U(2)}\times \frac{\SO(5)}{\U(2)}$:
on the left when the  symplectic form is not deformed, on the right after a deformation.}
\label{stellette}
\end{figure} 
}
\end{example}

\end{document}